\documentclass[graybox]{svmult}
\usepackage{mathptmx}       
\usepackage{helvet}         
\usepackage{courier}        
\usepackage{type1cm}        
\usepackage[varvw]{newtxmath}  
\usepackage{makeidx}         
\usepackage{graphicx}        
\usepackage{multicol}        
\usepackage[bottom]{footmisc}
\usepackage{hyphenat}
\makeindex             
\begin{document}

\title{Incorporating Shape Knowledge into Regression Models}
\author{Miltiadis Poursanidis, Patrick Link, Jochen Schmid and Uwe Teicher}
\institute{Miltiadis Poursanidis \at Fraunhofer Institute for Industrial Mathematics ITWM, 67663, Kaiserslautern, Germany, \email{miltiadis.poursanidis@itwm.fraunhofer.de}
\and Patrick Link \at Fraunhofer Institute for Machine Tools and Forming Technology IWU, 09126, Chemnitz, Germany, \email{patrick.link@iwu.fraunhofer.de}
\and Jochen Schmid \at Fraunhofer Institute for Industrial Mathematics ITWM, 67663, Kaiserslautern, Germany, \email{jochen.schmid@itwm.fraunhofer.de}
\and Uwe Teicher \at Fraunhofer Institute for Machine Tools and Forming Technology IWU, 09126, Chemnitz, Germany, \email{uwe.teicher@iwu.fraunhofer.de}}%
%
\maketitle

\abstract{Informed learning is an emerging field in machine learning that aims to compensate for insufficient data with prior knowledge. Shape knowledge covers many types of prior knowledge concerning the relationship of a function's output with respect to input variables, for example, monotonicity, convexity, etc. This shape knowledge -- when formalized into algebraic inequalities (shape constraints) --- can then be incorporated into the training of regression models via a constraint problem formulation. The defined shape-constrained regression problem is, mathematically speaking, a semi-infinite program (SIP). Although off-the-shelf algorithms can be used at this point to solve the SIP, we recommend an adaptive feasible-point algorithm that guarantees optimality up to arbitrary precision and strict fulfillment of the shape constraints. We apply this semi-infinite approach for shape-constrained regression (SIASCOR) to three application examples from manufacturing and one artificial example. One application example has not been considered in a shape-constrained regression setting before, so we used a methodology (ISI) to capture the shape knowledge and define corresponding shape constraints. Finally, we compare the SIASCOR method with a purely data-driven automated machine learning method (AutoML) and another approach for shape-constrained regression (SIAMOR) that uses a different solution algorithm.}

\section{Introduction}

Despite the success of machine learning (ML), purely data-driven machine learning models show limited performance when dealing with insufficient data. \index{small data} This is especially problematic in scientific and engineering contexts, where simulation or experimental data are costly in both time and resources \cite{Weichert.2019}. When the data set is small, machine learning models have difficulties to provide reliable models. The main issue is that the models do not behave as expected in regions with sparse or no data. When noise comes into play, this effect is even more severe as the models tend to learn spurious patterns from the data. In addition to that, it is difficult to measure the model's performance at sparse data regions. Also methods like cross-validation are often misleading because only the data set is considered for measuring the performance.

In many machine learning tasks, by contrast, there is additional prior knowledge \index{prior knowledge} available. Informed learning \index{Informed learning} emerged from the need to compensate for shortcomings in the data with supplementary prior knowledge \cite{Rueden.2021}. Machine learning models benefit from prior knowledge in various ways, but we highlight two in particular: interpretability \index{interpretability} and generalization. \index{generalization} Informed machine learning models are interpretable because they behave according to the imposed prior knowledge. For instance in production, there is usually high costs involved with wrong decision-making. Therefore, practitioners rely more on their knowledge than on the data, especially when data is sparse and noisy. For this reason, trustworthy prediction models \index{trustworthy machine learning} should incorporate prior knowledge to increase acceptance among practitioners. In science, interpretable models are key to the accumulation of scientific knowledge. In contrast to black-box models, theories can be developed based on these interpretable models that have known properties \cite{Karpatne.2017}. The other aspect is generalization, that is, the model's ability to achieve low errors on new data. Informed learning is expected to lead to improved generalization. Imposing prior knowledge gives control over regions in the domain with sparse data and makes models less prone to unexpected behavior. This typically leads to models that generalize better outside the data set. Another aspect is extrapolation. The authors in \cite{haider23} show that, in certain cases, shape constraints can lead to an improvement of the out-of-domain error. However, we do not consider extrapolation here.

In this work, we focus on prior knowledge concerning the qualitative shape of the model function. The definition of shape knowledge is very general and captures many properties such as boundedness or monotonicity of the model function. In the terminology of \cite{Rueden.2021}, shape knowledge can be categorized as either scientific knowledge (given in explicit formulas) or expert knowledge (common knowledge within a scientific field). Such shape knowledge can often be formulated as algebraic inequalities, so-called shape constraints.

First, we recapitulate the semi-infinite approach to shape-constrained regression (SIASCOR) \index{SIASCOR} from \cite{Link.2022}. Shape-constrained regression is a constraint problem formulation of a regression problem with the aim to incorporate shape constraints into the model function. Mathematically, this results in a so-called semi-infinite program (SIP) and can be solved, for instance, with adaptive feasible-point methods. In \cite{Link.2022}, the authors used the core algorithm from \cite{Schmid.2021} as their adaptive feasible-point algorithm to compute an approximate solution to the resulting SIP. In the present paper, by contrast, we use the simultaneous algorithm from \cite{Schmid.2021} as our adaptive feasible-point algorithm. In contrast to the core algorithm, it computes an approximate solution to the SIP of an arbitrary user-specified precision. According to the taxonomy of \cite{Rueden.2021}, SIASCOR integrates shape knowledge -- represented as algebraic inequalities -- into the training of the regression model. Second, we reconsider the methodology from \cite{Link.2022}. This methodology helps practitioners to capture shape knowledge in cooperation with experts in the field, to define shape constraints, and to incorporate the shape constraints into a regression model. We refer to this methodology as ISI from now on according to its three steps, namely the inspection, the specification, and the integration step. In terms of the taxonomy of \cite{Rueden.2021}, this can be partly viewed as a method that transforms expert knowledge into algebraic inequalities. In addition to that, the ISI methodology uses some shape-constrained regression method, for instance SIASCOR, to integrate these algebraic inequalities into the machine learning model.

We consider three real-world application examples from quality prediction in manufacturing: brushing, press hardening, and milling. On the one hand, all three examples have small and noisy data sets but, on the other hand, they can benefit from shape knowledge provided by experts. The brushing example was already considered in \cite{Link.2022} and the press hardening example in \cite{Kurnatowski.2021}, thus we reuse the same shape constraints as in these references. The milling case has not been studied in a shape-constrained regression setting before. Therefore, we use the ISI method to ensure that all shape knowledge is captured and, if possible, transformed into shape constraints. After that, we can apply SIASCOR with the obtained shape constraints. 
Moreover, and in contrast to \cite{Link.2022}, we compare SIASCOR to more sophisticated machine learning methods. We compare it with an automated machine learning method (AutoML) and  another semi-infinite approach to shape-constrained regression but with a different solution algorithm (SIAMOR) \cite{Kurnatowski.2021}. Another difference to \cite{Link.2022} is that here we use different settings of SIASCOR such as another solving algorithm of the SIP and anisotropic polynomial regression functions. We compare the resulting models in terms of shape compliance, training time, and cross-validated test error.

As an extension to the three real-world application examples, we introduce an artificial example to examine the generalization error. The generalization error is the error a model has on data not contained in the training set. Since our real-world examples have small data sets, we can not analyze the generalization error appropriately, especially in scarce data regions. We compare SIASCOR, for the first time, with AutoML, SIAMOR and with Ridge regression in terms of generalization error.

We organize the article as follows. In Section~\ref{sec:relatedwork}, we give a basic overview of the related work. In Section~\ref{sec:methods}, we introduce the informed machine learning approach SIASCOR for shape-constrained regression and present the methodology ISI to capture and integrate expert knowledge. Section~\ref{sec:realworldexamples} describes our three application examples from manufacturing -- namely press hardening, brushing, and milling -- and discusses the results of the comparative study. Section~\ref{sec:artificialexample} presents our artificial application example along with the analysis of the generalization error. Finally, we give a conclusion and an outlook in Section~\ref{sec:conclusion}.

\section{Related Work}
\label{sec:relatedwork}
In this section, we present some related work on informed learning, shape-constrained regression, semi-infinite programming, and expert-knowledge-based quality prediction in manufacturing.

Informed learning describes all approaches that incorporate prior knowledge into machine learning models. An overview of the field can be found in \cite{Rueden.2021}. In the present work, we focus on prior knowledge in the form of algebraic inequalities. According to the taxonomy of \cite{Rueden.2021}, algebraic inequalities are included in the class of algebraic equations. Algebraic inequalities can be integrated into regression models in four ways: When generating training data \cite{Ladicky.2015}, by restricting the hypothesis space \cite{Lu.2017, Bauckhage.2018}, during the learning algorithm \cite{Daw.31.10.2017, Stewart.19.09.2016, Diligenti.2017, Heese.2019, Kurnatowski.2021, Muralidhar.122018} and by modifying the final model \cite{Schmid.2020, Muralidhar.122018, King.2018}. Informed learning methods that integrate algebraic inequalities during training either treat the inequalities as soft constraints by adding a penalty term to the loss function \cite{Daw.31.10.2017, Stewart.19.09.2016, Diligenti.2017, Heese.2019} or as hard constraints by adding them as constraints to the loss minimization problem \cite{Kurnatowski.2021, Muralidhar.122018}. Among the different constraint types, shape constraints restrict the qualitative shape of the prediction function \cite{Gupta2020MultidimensionalSC}. One of the most prominent shape constraints is monotonicity and, in the literature, there already exist numerous approaches to enforcing monotonicity constraints during training \cite{Gupta.2016, Altendorf.2012, Lauer.2008.Regression, Chuang.2020, Riihimäki.2010}. Furthermore, the authors of \cite{AubinFrankowski.05.01.2021} consider various shape constraints in a kernel regression setting. They enforce the shape constraints on a finite set of points but sufficiently tighten the problem to fulfill the constraints on the entire input space. Polynomial shape-constrained regression is considered in \cite{Hall2018}, where the authors use SDP relaxations to solve the shape-constrained regression problem. The authors of \cite{Cozad.2015} and of \cite{Kurnatowski.2021} approach shape-constrained regression via semi-infinite programming. However, due to new mathematical results from the SIP community \cite{HatimDjelassi.}, more suitable algorithms can be used.

SIPs are optimization problems that have a finite number of decision variables and an infinite number of constraints. For an overview of the theory and how to handle the infinite constraints numerically, we refer to \cite{RembertReemtsen.1998, Hettich.1993, HatimDjelassi.}. Popular methods for solving SIPs are discretization methods \cite{Lopez.2007, Blankenship.1976} with the attention shifting to adaptive discretization. Among these, there is also a line of work concerning adaptive feasible-point methods \cite{Tsoukalas.2011, Mitsos.2011}, which guarantee termination at a feasible point. In \cite{Schmid.2021}, the authors leveraged the convexity -- a property inherent in most shape-constrained regression problems -- to obtain stronger results such as arbitrary optimality precision under weaker assumptions.

In \index{manufacturing} manufacturing, quality prediction \index{quality prediction} is used to both monitor product quality and optimize processes. Quality prediction models are either data-driven or rely on physical equations. In the context of manufacturing, the use of data-driven models is a challenging task due to data scarcity. As mentioned in \cite{Weichert.2019}, complex models are applied to describe complex relationships that have few data available. The problem of data scarcity is also reported in other domains, such as process engineering \cite{Napoli.2011}.
In order to handle small data sets, multi-model approaches \cite{Li.2012, Chang.2015} or polynomial chaos expansion \cite{Torre.2019} have been used. Another technique is to generate artificial data via bootstrapping \cite{Napoli.2011, Tsai.2008} or mega-trend-diffusion \cite{Li.2007, Li.2013}. Besides these data-driven methods, there are also expert-knowledge-based approaches \cite{Hall2018}. In manufacturing, the most common sources of knowledge are scientific and expert knowledge, according to \cite{Rueden.2021}. The authors of \cite{Zhang.2020, Lokrantz.2018, He.2019} integrate probabilistic relationships into the hypothesis sets of Bayesian networks. In addition,  \cite{Lu.2017} and  \cite{Nagarajan.2019} restrict the hypothesis set of neural networks with algebraic equations and  knowledge graphs, respectively. The authors of \cite{Hao.2020} incorporate algebraic equations into the training of Gaussian process models.

\section{Methods}
\label{sec:methods}
In the first subsection, we describe a methodology (SIASCOR) that integrates shape constraints into regression models via semi-infinite programming. In practice, however, shape knowledge is not always available as algebraic inequalities. Usually, there is merely the expert's intuition that needs to be captured and formalized into shape constraints. Therefore, in the second subsection, we recall a methodology (ISI) that helps to capture shape knowledge \index{shape knowldge} and convert it into algebraic inequalities. These algebraic inequalities can then be integrated into the regression model using SIASCOR, for instance.

\subsection{SIASCOR}
\label{sec:SIASCOR}
The goal of classical regression settings usually boils down to finding a model function that fits some data.  Assume we have additional prior knowledge about the shape of the  input-output relationship to be learned. When shape knowledge is formalized in terms of inequality constraints, we call them shape constraints. \index{shape constraints}  Common forms of shape knowledge, for instance, is monotonicity or convexity. The corresponding shape constraints restrict the first or second partial derivative with respect to some input variable of the model function to be positive. Then, the goal of shape-constrained regression \index{shape-constrained regression} is to find a model function that both fits the data and complies with the shape constraints.

Suppose we are given some data set $\mathcal D = \{ (x^k, y^k) \in X \times \mathbb R : k = 1, \dots , n\}$ consisting of input data points $x^k \in \mathbb R^d$ and output data points $y^k \in \mathbb R$. We assume that the (unknown) input-output relationship to be learned can be represented by model functions $\widehat{y}_w: X \to \mathbb R$ of the form $\widehat{y}_w(x) := w^{\top} \phi(x)$. In other words, we take our hypothesis space to be the set of functions $\widehat{y}_w$ with $w \in W$. In the above formula, $x \in X$ and $w \in W$ denote the input variables and model  parameters, respectively, and we assume the input-variable and model-parameter spaces $X \subset \mathbb R^d$ and $W \subset \mathbb R^m$ to be compact, convex sets. Also, $\phi: X \to \mathbb R^m$ is a feature mapping that is sufficiently often differentiable. We further assume that the shape constraints can be expressed by constraining a function $g_i$ that is given in terms of affine-linear combinations of the partial derivatives of the model function $\widehat{y}_w$. Clearly, boundedness, monotonicity or convexity constraints \index{monotonicity constraints} \index{convexity constraint} can be cast in this form; for instance, monotonic increasingness w.r.t. $x_j$ is equivalent to the condition $g(w,x) :=  - w^T \partial_{x_j} \phi (x)\leq 0$. Consequently, in mathematical terms, shape-constrained regression problems take the form 
\begin{equation}
	\label{eq:shape-constrained-regression}
	\begin{aligned}
		\min_{w\in W} &\sum_{k=1}^n |y^k - \widehat{y}_w(x^k)|^2 + \lambda ||w||^2 \\
		\text{s.t.} \ &g_i(w, x)\leq 0 \ \text{for all } x\in X \text{ and } i \in I,
	\end{aligned}
\end{equation} 
where $\lambda > 0$ is some regularization parameter, $\|w\|^2$ denotes the squared $\ell^2$-norm of the model parameter $w$ (ridge regularization) and $I$ is a finite index set indexing the shape constraints. Note that we restricted ourselves here to the ridge regression and to model functions that are linear w.r.t. to their model parameters. For more general cases, see \cite{Schmid.2021}.

Problem \eqref{eq:shape-constrained-regression} is a so-called semi-infinite program (SIP) \index{semi-infinite program}. Assume the problem is feasible, i.e. there exists a $w\in W$ such that $g_i(w, x) \leq 0$ for all $i\in I$ and $x\in X$. Intuitively, this means that there exists a function from our hypothesis space that satisfies all shape constraints. Then, problem \eqref{eq:shape-constrained-regression} has a unique solution, by its strict convexity and our continuity assumptions on $\phi$. See \cite{Schmid.2021, Hettich.1993}, for example. There exist many approaches for solving SIPs \cite{HatimDjelassi.} and in particular convex SIPs \cite{RembertReemtsen.1998}. We prefer feasible-point methods because they guarantee termination at a feasible point. Other methods mostly guarantee feasibility only as the iteration number tends to infinity. Hence, we use the simultaneous algorithm from \cite{Schmid.2021}, a feasible-point method that leverages the convexity of the problem to provide guarantees for approximate, feasible solutions, while the only assumption is strict feasibility of problem \eqref{eq:shape-constrained-regression}. There are more algorithmic merits of the approach but we will not detail them here. Note that the algorithm we used in this paper is different from \cite{Link.2022} where we used the core algorithm from \cite{Schmid.2021}. The core algorithm also guarantees feasibility but does not guarantee optimality of arbitrary precision. Besides, it is different from \cite{Cozad.2015} where the authors do not use a feasible-point method in the first place.

After having defined the method more precisely, we can see how it can be embedded into the taxonomy of \cite{Rueden.2021}: knowledge, in our context, is given in the form of shape constraints which, ultimately, are algebraic inequalities. Then, these algebraic inequalities are integrated during the learning algorithm through a constrained optimization problem formulation.

In contrast to soft-constraint methods, SIASCOR imposes hard constraints \index{hard constraints} on the model function. This is suitable for applications where the model function needs to satisfy the constraints strictly, for instance when the model needs to be in accordance with physical laws. Despite that, one can relax the constraints by a small value  $\varepsilon > 0$ if the constraints do not need to be fulfilled strictly. This can be done by subtracting the value  $\varepsilon$ from all shape constraint functions $g_i$. In this case, the resulting model function complies with the relaxed shape constraints.

\subsection{ISI} 
\label{sec:ISI}

In the previous section, we described how SIASCOR incorporates shape constraints into the training of machine learning models of the form $\widehat{y}(x) = w^T\phi(x)$. In practice, these shape constraints need to be developed in collaboration with experts in the field. The ISI \index{ISI} methodology supports practitioners in capturing shape expert knowledge, in formalizing it into shape constraints and in producing a shape-compliant model. In this section, we summarize the ISI methodology that was introduced in \cite{Link.2022}.

\begin{figure}[t]
	\begin{center}
		\makebox[\textwidth]{\includegraphics[width=0.8\textwidth]{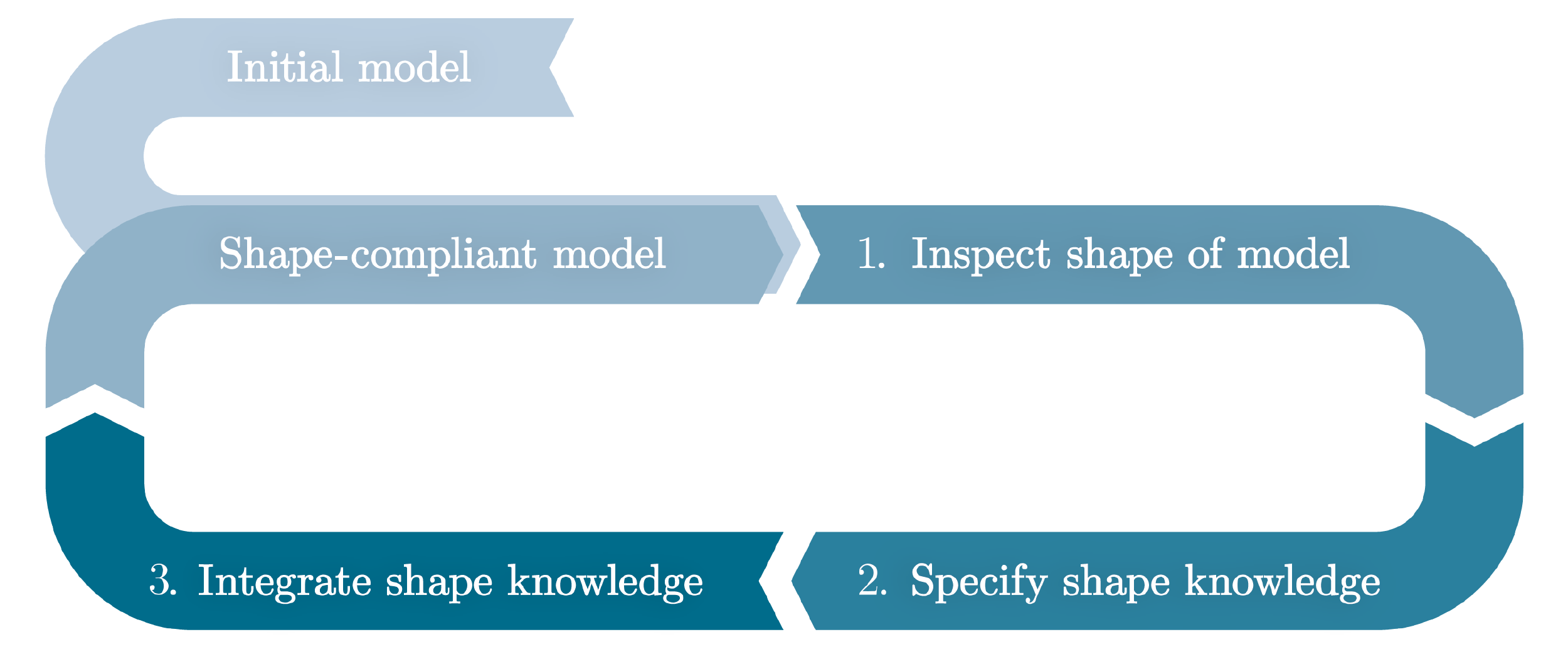}}
	\end{center}
	\caption[Schematic]%
	{Schematic of the three-step methodology  with an initial purely data-based model as its input and a shape-knowledge-compliant model as its output. Image source: \cite{Link.2022}.}
	\label{fig:methodology-schematic} 
\end{figure}

The schematic procedure of ISI is depicted in Figure~\ref{fig:methodology-schematic}. As its input, the methodology requires an initial model $\widehat{y}^0$ that may be purely data-based. Then, the methodology proceeds in the following three steps:

\begin{enumerate}
	\item Inspection of the initial model by a process expert
	\item Specification of shape expert knowledge by the expert
	\item Integration of the specified shape expert knowledge into the training of a new model
\end{enumerate}

These three steps generate a shape-compliant model. If this model does not behave as the expert expects, the three steps can be repeated. Note that in the second iteration one must take the shape-compliant model from the first iteration as the input for the inspection. This can be repeated as often as necessary while always taking the current shape-constrained model as the input model for the next iteration. 

The initial model is the starting point for the introduced methodology. The initial model visually supports the expert in analyzing the relationship between the inputs and outputs. In principle, any type of model function can be used as an initial model. However, we recommend choosing the initial model from the same hypothesis space as the SIASCOR model. In other words, the initial model should be of the form $\widehat{y}^0 (x) = w_0^T\phi(x)$ for some $w_0\in W$. This way, the expert can get an intuition of how model functions from this hypothesis space behave when they are not shape-constrained, especially in regions with sparse data. The parameter $w_0\in W$ can be computed in any purely data-based regression method, such as ridge regression.

After generating an initial model, the first step of the methodology is to provide the expert with one- and two-dimensional graphs at multiple points of high- and low-fidelity. Here, a custom notion of fidelity can be used. For instance, one can consider points $x_\text{high}, x_\text{low}\in X$ with maximal and minimal distance to all data points, respectively. These one- and two-dimensional graphs at the high- and low-fidelity points help the expert understand the functional relationships behind the data by exploring the initial model along the input space. This way, the expert can find shape behavior that either confirms or contradicts their intuition.

\begin{figure}[t]
	\label{fig:shape_constraints}
	\sidecaption[t]
	\includegraphics[scale=.7]{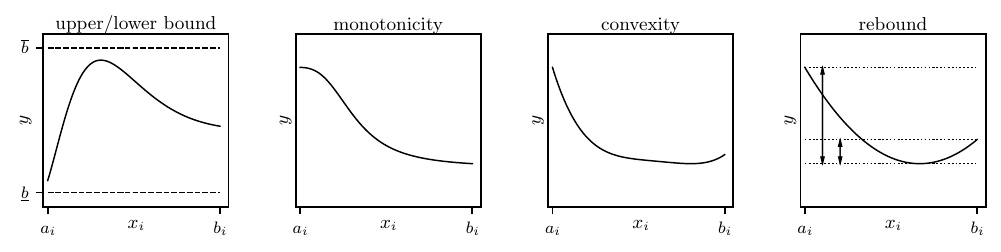}
	%
	%
	\caption{Sketches of the the shape constraints we considered in our applications. The upper- and lower-bound constraints bound the maximal and minimal values of the function, respectively. Monotonicity constraints (either increasing or decreasing)  restrict the first derivative of the function to be either positive or negative. Similarly, convexity constraints assume the second derivative to be positive (or negative if concave). The rebound constraint describes, loosely speaking, how much a function can increase again after it has witnessed a decrease. For a more details, see \cite{Link.2022} and for proofs \cite{Schmid.2021}.}
	\label{fig:expert-knowledge}       
\end{figure}

In the second step of the methodology, the process expert specifies the shape expert knowledge, that is, his intuition about the qualitative functional relationship of the output along different dimensions of the input space. Then, shape knowledge is converted into shape constraints. The expert can choose from a variety of common shape knowledge with associated shape constraints such as monotonicity, convexity or concavity or upper and lower bounds, see Figure \ref{fig:expert-knowledge} for a selection of qualitative shape knowledge. Also, multiple shape constraints can be combined. In case that the shape knowledge cannot be composed by common shape constraints, practitioners may consider designing a new shape constraint. This, for instance, was the case in \cite{Link.2022} that resulted in defining the rebound constraint.

With the shape constraints at hand, the third step of the procedure is to find a shape-compliant model. At this point, SIASCOR can be used but can be interchanged with any shape-constrained regression model. In comparison to the initial model, the shape-compliant model fits both the data and satisfies the defined shape constraints.

After the first iteration, the ISI methodology produces a model that is compliant with the shape constraints imposed in the first iteration. Yet, we can neither guarantee that the imposed shape constraints are complete nor that they are all correct. It can happen, for instance, that some shape knowledge has been overlooked in the first run and, therefore, necessary shape constraints have not been imposed. We witnessed this when we applied the ISI methodology in a brushing use case in \cite{Link.2022}. Similarly, it can happen that some of the proposed shape constrains are not in harmony with the data. We witnessed this case in the milling application example that we detailed in Section \ref{sec:milling}. Therefore, we repeat the ISI procedure by replacing the initial model with the refined shape-compliant model as often as needed until the expert detects no more inconsistencies between the shape of the model function and his shape knowledge.

\section{Real-world Application Examples}
\label{sec:realworldexamples}

In this section, we present three real-world application examples for shape\hyp{}constrained regression. We apply the SIASCOR method from Section \ref{sec:SIASCOR} to incorporate prior knowledge given in the form of shape constraints. Fpr two application examples, the shape constraints have been already captured in previous works, using the ISI methodology -- either implicitly in the case of \cite{Kurnatowski.2021} or explicitly in the case of \cite{Link.2022}. The other example has not been considered in a shape-constrained regression setting before, so we did not have any shape constraints available. Therefore, we used the ISI methodology to capture the prior knowledge and to define the shape constraints, together with experts in the field. Finally, we compare SIASCOR with both a data-driven machine learning model and another shape-constrained regression approach (SIAMOR) \cite{Kurnatowski.2021}. 

Every subsection deals with a distinct application example from quality prediction for manufacturing processes. We consider three manufacturing processes: press hardening, brushing, and milling. Quality prediction for the press hardening case has already been considered in \cite{Kurnatowski.2021} and for the brushing case in \cite{Link.2022}. In these two cases, the shape constraints have already been identified and we can directly apply SIASCOR with respect to the shape constraints proposed by the experts. In contrast, the milling case has not been considered for shape-constrained regression so far. Hence, we first apply the ISI procedure in collaboration with experts in the field to capture their shape knowledge and formalize it in terms of shape constraints. With the shape constraints at hand, we finally apply SIASCOR. 

First, we compare SIASCOR with a purely data-driven machine learning approach to highlight the necessity of imposing shape constraints. Concretely, we compare it with auto-sklearn \cite{Feurer.2015}, an sklearn \cite{Pedregosa.2011} package for automatic machine learning (AutoML). We choose an AutoML approach, because this has become a preferred choice in many fields of machine learning \cite{He.2021} due to its systematic and automated manner of model and hyperparameter selection. We consider two versions of auto-sklearn in our comparison: one with default settings and another with handpicked regression models, data- and feature preprocessing to avoid unphysical behavior through discontinuities in the shape of the model function. More precisely, we exclude ``extra\_trees\_preproc\_for\_regression'' and ``random\_trees\_embedding'' from the feature preprocessing and ``adaboost'', ``decision\_tree'', ``extra\_trees'', ``gradient\_boosting'',  ``k\_nearest\_neighbors'' and ``random\_forest'' from the regressor selection.

Second, we compare SIASCOR with its predecessor SIAMOR -- an alternative method for shape-constrained regression. It was initially developed to enforce monotonicity of regression models, but it can be easily adapted to general shape constraints. Both methods consider a SIP formulation of shape-constrained regression. In contrast to SIASCOR, however, SIAMOR uses a different algorithm to solve the SIP which does not provide a theoretical guarantee for feasibility. It merely guarantees feasibility on a reference grid and, generally, only after infinitely many iterations. In other words, we can only guarantee shape compliance on a grid of finitely many points in the limit of infinitely many iterations. In all applications, we chose a reference grid of 20 points per dimension. Besides the theoretical differences, we compare these two methods in practice.

We consider three criteria in our comparison: shape compliance, training time and test error. For all three criteria, we conducted a 10-fold cross-validation. First, we counted how many shape constraints the final models violated in each fold and then we averaged these numbers over all folds. For every shape constraint we sampled 10000 random points in $X$ and tested the shape compliance on 100 equidistant points along the relevant axis. Besides, we visualize some of the violations at selected points. During the training of each fold, we also measured the wall clock time and averaged it across all folds. The training was conducted on a standard office laptop. Lastly, we measured the root-mean-squared errors RMSE of the model on all test sets and, again, averaged along all folds.

\subsection{Press Hardening}
\label{sec:presshardening}

\begin{table}[!t]
	\label{tab:press_hardening}
	\caption{Press hardening case: Comparison of SIASCOR, Auto-Sklearn with default settings, Auto-Sklearn with custom settings, and SIAMOR. The table lists the average number of shape constraints that the models violate after trained on each fold of a 10-fold cross-validation. We considered three shape constraints in total. Moreover, both training time and test error were measured in each fold and then averaged over the conducted cross-validation. In addition to the mean of the cross-validated (CV) test errors we added also the standard deviation. The training time is given in h:m:s and the error in root mean squared error (RMSE).} 
	\label{tab:press_hardening}       
	\centering
	\begin{tabular}{p{3cm}p{2.75cm}p{2.75cm}p{2.75cm}}
		\hline\noalign{\smallskip}
		Model & CV Test Error & Training Time & Shape Violations \\
		\noalign{\smallskip}\svhline\noalign{\smallskip}
		SIASCOR &  17.92 $\pm$ 5.92 & 00:03:09 & 0 out of 3\\
		AutoML1 & 16.06 $\pm$ 6.26 & 00:09:56 & 3 out of 3\\
		AutoML2 & 15.73 $\pm$ 5.84 & 00:09:56 & 3 out of 3\\
		SIAMOR  & 17.92 $\pm$ 5.92 & 00:02:14 & 1 out of 3\\
		\noalign{\smallskip}\hline\noalign{\smallskip}
	\end{tabular}
\end{table}

Our first application example is press hardening, a hot sheet metal forming process \cite{Neugebauer.2012}. During the forming process, the hot sheet metal is formed and subsequently quenched to achieve improved hardness of the parts. The goal in this application is to build a model that predicts the hardness of the metal sheet as a function of four process parameters, namely the furnace temperature $T_f$, the handling time from furnace to press $\Delta t_h$, the press force $F_p$ and the quenching time $t_q$. For the informed learning task, we use the 60 data samples that were generated in an experimental setup and the defined monotonicity constraints, both from \cite{Kurnatowski.2021}.

Now we present our parameter settings for SIASCOR. First, we scaled the input data to the unit cube based on the ranges provided in \cite{Kurnatowski.2021}. This way, we can consider $X = [0, 1]^4$ and all the data is contained in that domain $\mathcal{D}\subset X \times \mathbb R$. Moreover,  we choose $\phi$ to be an anisotropic polynomial feature map. This means that the maximal polynomial degree can be different for each of the four dimensions. Specifically, we chose the maximal degrees to be (4, 1, 3, 4) for $(T_f, \Delta t_h, F_p, t_q)$, as suggested by the process experts. Accordingly, we set the parameter space to be $W=[-10^5, 10^5]^{54}$ and, furthermore, the regularization parameter  $\lambda = 0.0001$. The shape constraints from \cite{Kurnatowski.2021} are monotonic increasingness in $T_f$ and $t_q$ and monotonic decreasingness in $\Delta t_h$. The settings above specify the SIP in \eqref{eq:shape-constrained-regression}. Afterwards, we apply the main algorithm from \cite{Schmid.2021} to solve the SIP with optimality precision $\delta=0.0001$. We are not going into detail here which parameters values we chose for the adaptive feasible-point algorithm. For more details on the algorithm and its parameter settings, see \cite{Schmid.2021}.

For both AutoML approaches, we set the strategy that chooses the best model to be a 10-fold cross-validation and fixed the maximum search time to 600 seconds, as suggested in \cite{Feurer.08.07.2020}. Afterwards, the AutoML models were refit on the entire training set. For SIAMOR, we set everything as in \cite{Kurnatowski.2021}, especially the size of the reference grid being 20 points per dimension. 

\begin{figure}[t]
	\sidecaption[t]
	\includegraphics[scale=.375]{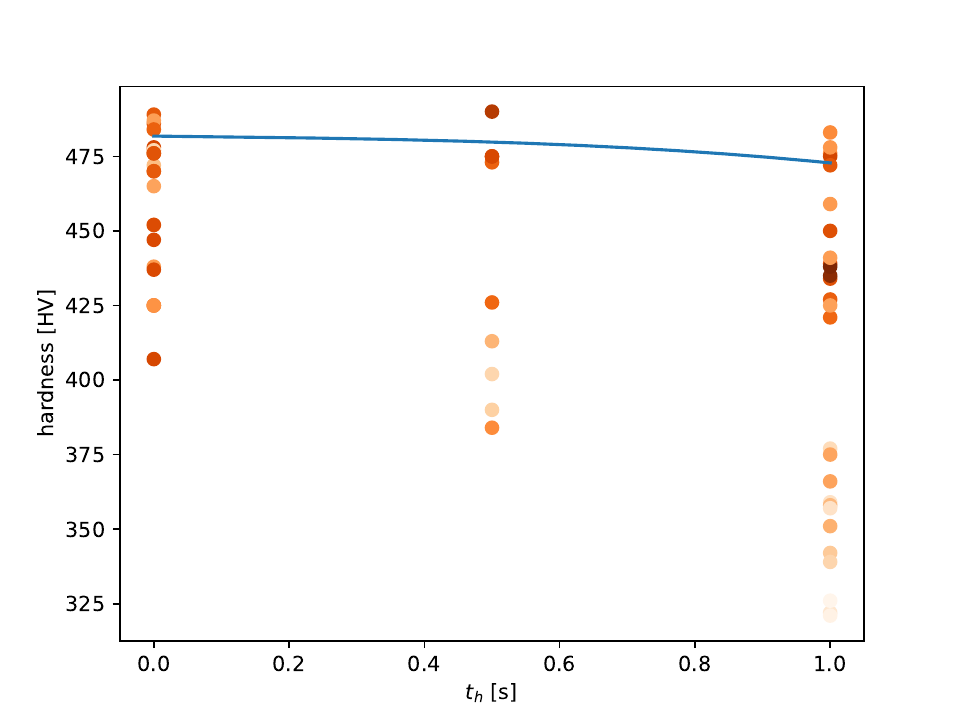}
	\includegraphics[scale=.375]{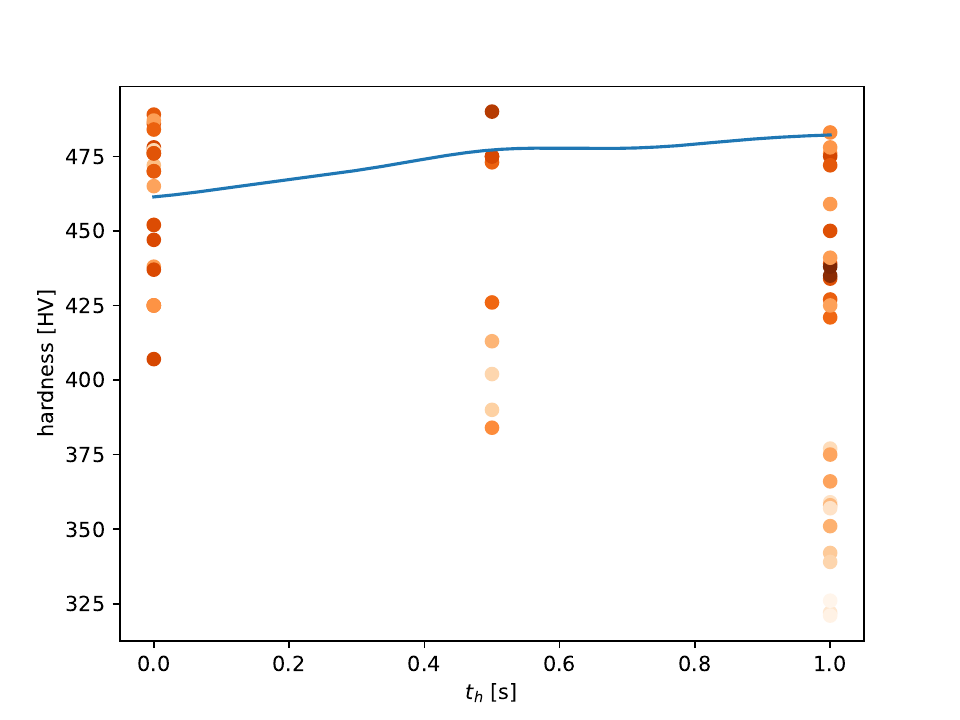}
	%
	%
	\caption{The left  graph depicts the shape of the SIASCOR model and the right graph the shape of the auto-sklearn model both along the $t_c$ variable fixed at the point $x=(0.998, 0.506, 0.154, 0.334)$.  The black points represent data points located along the axis and the orange points are the remaining data points projected onto the $\Delta t_h$ dimension. The darker the orange points are, the closer the Euclidean distance is between the data point and its projection. Both models were trained on 90\% of the data that are visible in the graph.}
	\label{fig:violation-press-hardening}       
\end{figure}

Table \ref{tab:press_hardening} shows that both AutoML versions produce models that are not shape-compliant. In fact, we observe that all models produced during the cross-validation violate all three imposed shape constraints. Also the trained SIAMOR models violate one out of four shape constraints on average. As expected, the SIASCOR model satisfies all shape constraints. Figure \ref{fig:violation-press-hardening} juxtaposes one-dimensional graphs of the SIASCOR and the auto-sklearn model at the same point. We can see graphically that the SIASCOR model is monotonically decreasing while the auto-sklearn model is not. However, training the SIASCOR model takes a little bit longer than the other models, which can be traced back to the computationally expensive global optimization steps that guarantee shape compliance. Furthermore, we infer from Table \ref{tab:press_hardening} that the averaged error on the test sets during a cross-validation is more or less the same, considering the scale of the hardness (roughly between 300-500). The AutoML model has slightly lower test errors but considering the amount of data, this does not imply good generalization outside the data.

\subsection{Brushing}
\label{sec:brushing}

The brushing process is a metal-cutting process used for machining of surface structures with the help of brushes. In quality prediction, the goal is to predict the surface roughness given adjustable process parameters. In our example, we had a data set consisting of 125 points which were generated in an experimental setup where the average arithmetic roughness $R_a$ was measured for various settings of five process parameters: the diameter of the abrasive grits $Dia$, the cutting time $t_c$, number of revolutions of the brush $n_b$, number of revolutions of the work piece $n_w$, and the cutting depth $a_e$. Aside from the data, experts have prior knowledge about the behavior of these machining processes. In \cite{Link.2022}, we used the ISI method to formalize shape constraints from the expert knowledge. The shape constraints that the experts suggested are visualized Figure \ref{fig:expert-knowledge}. In short, we had upper and lower bounds, monotonicity, convexity, and rebound constraints. For a mathematical description of the shape constraints, see \cite{Link.2022}.

Similar to the press hardening case, we scale the input data to the unit cube so that we can consider  $X = [0, 1]^5$ and set the regularization parameter $\lambda = 0.01$. We choose $\phi$ to be an anisotropic polynomial feature map with degrees (1, 5, 2, 2, 2) for $(Dia, t_c, n_b, n_w, a_e)$ and the parameter space to be $W=[-10^5, 10^5]^{136}$. Additionally, we integrated the same shape constraints as in \cite{Link.2022}. With these settings, the SIP in \eqref{eq:shape-constrained-regression} is specified for both SIAMOR and SIASCOR.  Afterwards, we apply the algorithm from \cite{Schmid.2021} to solve the SIP and set the optimality precision to $\delta=0.0001$ and SIAMOR with the same settings as in the press hardening case.

\begin{table}[!t]
	\caption{Brushing case: Comparison of SIASCOR, Auto-Sklearn with default settings, Auto-Sklearn with custom settings, and SIAMOR. The table lists the average number of shape constraints that the models violate after trained on each fold a 10-fold cross-validation. We considered ten shape constraints in total. Moreover, both training time and test error were measured in each fold and then averaged over the conducted cross-validation. In addition to the mean of the cross-validated (CV) test errors we added also the standard deviation. The training time is given in h:m:s and the error in root mean squared error (RMSE).}
	\label{tab:brushing}       
	%
	%
	\begin{tabular}{p{3cm}p{2.75cm}p{2.75cm}p{2.75cm}}
		\hline\noalign{\smallskip}
		Model & CV Test Error & Training Time & Shape Violations \\
		\noalign{\smallskip}\svhline\noalign{\smallskip}
		SIASCOR & 0.0272 $\pm$ 0.008 & 01:10:51 & 0 out of 10\\
		AutoML1 & 0.0216 $\pm$ 0.006 & 00:09:56 & 7.8 out of 10\\
		AutoML2 & 0.0217 $\pm$ 0.006 & 00:09:56 & 7.8 out of 10\\
		SIAMOR  & 0.0267 $\pm$ 0.008 & 00:44:15 & 1 out of 10\\
		\noalign{\smallskip}\hline\noalign{\smallskip}
	\end{tabular}
\end{table}

\begin{figure}[t]
	\sidecaption[t]
	\includegraphics[scale=.375]{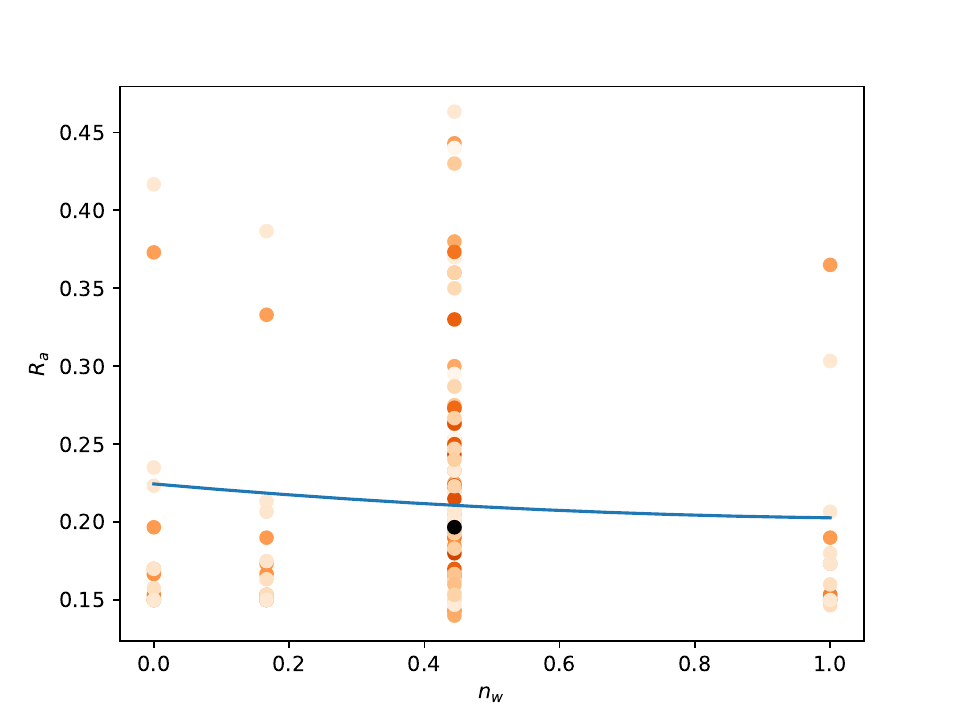}
	\includegraphics[scale=.375]{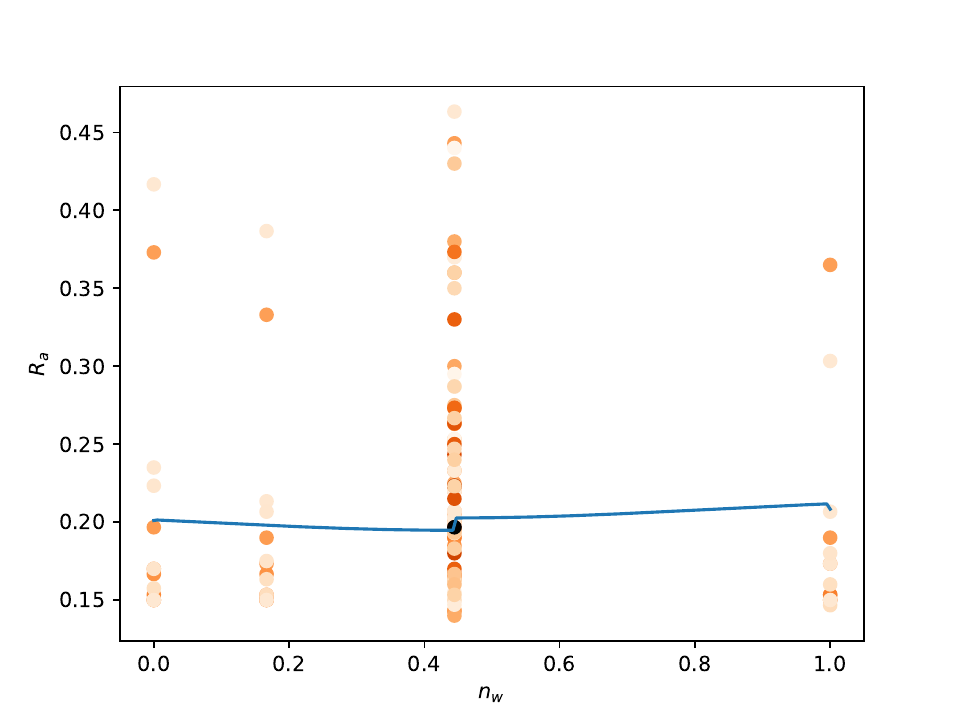}
	\caption{The left  graph depicts the shape of the SIASCOR model and the right graph the shape of the auto-sklearn model both along the $n_w$ variable fixed at the point $x=(1  ,0.5, 0.667 ,0.444, 0.334)$.  The black points represent the data points located along this axis and the orange points are the projected data points onto the $n_w$ dimension. The darker the orange points are, the closer the Euclidean distance is between the data point and its projection. Both models were trained on 90\% of the data that are visible in the graph.}
	\label{fig:violation-brushing}      
\end{figure}

Table \ref{tab:brushing} shows the results for the brushing case. Again, both AutoML approaches violate on average 7.8 out of ten shape constraints. Even the final model of SIAMOR violates on average 1 out of ten shape constraints. But, as expected, all SIASCOR models are in accordance with all shape constraints. In Figure \ref{fig:violation-brushing}, we see one exemplary shape violation of an AutoML model. On the right-hand side, the AutoML model violates the rebound and the convexity constraint and is, hence, not in accordance with physical laws. Although there are shape violations, the AutoML model used for Figure \ref{fig:violation-brushing} does not perform too badly. Similar to the convexity constraint in the graph, all other shape constraints were violated only in a very small region. This suggests an alternative measure for shape compliance that also takes into account the size of the region where violations occur. However, we do not go into that here. Anyway, this is also a representative example that AutoML may sometimes violate shape constraints but, nevertheless, not be entirely catastrophic. Moreover, the average training time of SIASCOR was higher compared to the other methods. The averaged test error is so low for all models considering that the data ranges from 0.14 to 0.46333.

\subsection{Milling}
\label{sec:milling}

Milling processes are characterized by high flexibility and productivity, which is why they are used for precision applications with high performance \cite{TEICHER201951}. Milling is a cutting process in which the tool rotates based on a geometrically determined definition and is subjected to chip removal from a workpiece due to the superposition of two effective directions of cutting and feed direction \cite{Einordnung}.
An assessment of the process is performed mainly on the basis of the quality with the help of roughness parameters (mainly $R_a$) and the mechanical loads on the basis of the parameter of the cutting force $F_c$. This addresses the decisive economic parameters for quality and energy consumption, so that their predictive capability is of high relevance. Technically, this is addressed by means of coating systems, which represent a central key to improving cutting properties due to the coating tribology \cite{PRENGEL1998183}. Aluminum represents a central research field for optimization due to its wide range of applications and the alloy-dependent variability of technical properties during milling \cite{Teichert19}.
Using the example of milling various aluminum alloys with coated solid carbide tools and varying the technological parameters of cutting speed $v_c$ and tooth feed $f_z$, the influence on the arithmetic center-line depth $R_a$ and cutting force $F_c$ is analyzed. In addition, the parameter of the friction coefficient $\mu$ of the coating system was taken into account, although this parameter is a one-dimensional quantity and does not fully reflect the tribological properties of the system. In total, we consider two outputs ($R_a$ and $F_c$) and three materials (``EN-AW5754'', ``EN-AW6082'', and ``EN-AW7075''), resulting in six different informed learning tasks. For each material and each output we had 80 data points.
This milling use case has not been considered for shape-constrained regression before. Thus, we had no shape constraints available in the beginning. So, we applied the ISI methodology in order to gather the shape constraints for every task. In the first iteration of the methodology, we trained a purely data-driven polynomial model on the data. After the inspecting step, the expert provided us with an initial set of shape constraints. Then, we incorporated all these shape constraints into the prediction model using the SIASCOR method. During the inspection step in the second ISI iteration, however, we noticed that the model had a different shape than expected along the $\mu$-axis. After taking a closer look at the data, we realized that the imposed shape constraints and the data were in conflict.  As said before, the one-dimensional variable $\mu$ does not fully capture the tribological properties and its influence on the variables $R_a$ and $F_c$ is hard to interpret. We therefore decided to drop the shape constraints along $\mu$ and left the formation of shape up to the data. In the inspection step of the following -- third -- ISI iteration, the expert was satisfied with the overall shape of the prediction function. And so, the iterative ISI procedure was finished at that point.

\begin{table}[!t]
	\caption{Milling case for output $F_{c}$: Comparison of SIASCOR, auto-sklearn with default settings, auto-sklearn with custom settings, and SIAMOR for all three materials. The table lists the average number of shape constraints that the models violate after trained on each fold a 10-fold cross-validation. We considered six to seven  shape constraints for the various materials. Moreover, both training time and test error were measured in each fold and then averaged over the conducted cross-validation. The training time is given in h:m:s and the error in root mean squared error (RMSE).}
	\label{tab:1}       
	%
	%
	\begin{tabular}{p{2.4cm}p{2.2cm}p{2.2cm}p{2.2cm}p{2.2cm}}
		\hline\noalign{\smallskip}
		Material & Model & CV Test Error & Training time & Shape Violations \\
		\noalign{\smallskip}\svhline\noalign{\smallskip}
		AW5754 & SIASCOR & 6.98 $\pm$ 3.52 & 00:00:19 & 0 out of 6\\
		AW5754 & AutoML1 & 6.99 $\pm$ 3.73 & 00:09:56 & 4 out of 6\\
		AW5754 & AutoML2 & 7.13 $\pm$ 3.87 & 00:09:56 & 4 out of 6\\
		AW5754 & SIAMOR  & 7.01 $\pm$ 3.50 & 00:00:09 & 0.2 out of 6\\
		\noalign{\smallskip}\hline\noalign{\smallskip}
		AW6082 & SIASCOR & \raggedright 7.86 $\pm$ 1.71 & 00:00:09 & 0 out of 7\\
		AW6082 & AutoML1 & 7.60 $\pm$ 4.16 & 00:09:56 & 5 out of 7\\
		AW6082 & AutoML2 & 7.57 $\pm$ 2.22 & 00:09:56 & 5 out of 7\\
		AW6082 & SIAMOR  & 7.90 $\pm$ 1.69 & 00:01:45 & 0.2 out of 7 \\
		\noalign{\smallskip}\hline\noalign{\smallskip}
		AW7075 & SIASCOR & 19.73 $\pm$ 11.08 & 00:00:12 & 0 out of 7\\
		AW7075 & AutoML1 & 20.27 $\pm$ 11.22 & 00:09:56 & 5 out of 7\\
		AW7075 & AutoML2 & 19.94 $\pm$ 11.40 & 00:09:56 & 5 out of 7\\
		AW7075 & SIAMOR & 19.74 $\pm$ 11.05 & 00:00:03 & 0 out of 7\\
		\noalign{\smallskip}\hline\noalign{\smallskip}
	\end{tabular}
\end{table}

Now we specify the parameter settings for SIASCOR used in the second and final iteration of ISI. First, we scale the input variables to the unit cube, according to the ranges $\mu \in [0.07, 0.5]$, $v_c\in [100, 1000]$, and $f_z\in [0.025, 0.25]$. The feature map is again an anisotropic polynomial with degrees (3, 2, 3) for ($\mu$, $v_c$, $f_z$) that induces the parameter space $W=[-10^5, 10^5]^{20}$.  In addition, we set the regularization parameter to $\lambda = 0.00001$ for all models. Now let us consider the shape constraints for the models with output  $F_c$. We impose a lower bound constraint with value 0 and an upper bound constraint with value 180,  a decreasingness constraint along $v_c$, and convexity constraints along every dimension. Moreover, the models for material EN-AW6082 and EN-AW7075 had an additional increasingness constraint along $f_z$. For the output $R_a$,  we impose a lower bound constraint with value 0, an upper bound constraint with value 6, and convexity constraints for all dimensions. Moreover, the models for material EN-AW5754 and EN-AW7075 had an increasingness constraint along $f_z$ and the model for material EN-AW7075 a decreasingness constraint along $v_c$. These settings specify the SIP in \eqref{eq:shape-constrained-regression}. Afterwards, we apply the main algorithm from \cite{Schmid.2021} to solve the SIP with optimality precision $\delta=0.0001$ for the models predicting $R_a$ and $\delta = 1$ for the ones predicting $F_c$. Furthermore, we used the same shape constraints as above for SIAMOR. Apart from that, the SIAMOR and the AutoML models were trained with same settings as in the other two application examples.

\begin{table}[t]
	\caption{Milling case for output $R_a$: Comparison of SIASCOR, auto-sklearn with default settings, auto-sklearn with custom settings, and SIAMOR for all three materials. The table lists the average number of shape constraints that the models violate after trained on each fold a 10-fold cross-validation. We considered five to seven  shape constraints for the various materials. Moreover, both training time and test error were measured in each fold and then averaged over the conducted cross-validation. In addition to the mean of the cross-validated (CV) test errors we added also the standard deviation. The training time is given in h:m:s and the error in root mean squared error (RMSE)}
	\begin{tabular}{p{2.4cm}p{2.2cm}p{2.2cm}p{2.2cm}p{2.2cm}}
		\hline\noalign{\smallskip}
		Material & Models & Test Error & Training time & Shape Violations \\
		\noalign{\smallskip}\svhline\noalign{\smallskip}
		AW5754 & SIASCOR & 0.3608 $\pm$ 0.1934 & 00:00:07 & 0 out of 6\\
		AW5754 & AutoML1 & 0.3163 $\pm$ 0.1498 & 00:09:58 & 4 out of 6\\
		AW5754 & AutoML2 & 0.3164 $\pm$ 0.1474 & 00:09:58 & 4 out of 6\\
		AW5754 & SIAMOR  & 0.3611 $\pm$ 0.1936 & 00:00:02  & 0 out of 6\\
		\noalign{\smallskip}\hline\noalign{\smallskip}
		AW6082 & SIASCOR & 0.2480 $\pm$ 0.0607 & 00:00:06 & 0 out of 5\\
		AW6082 & AutoML1 & 0.2292 $\pm$ 0.0859 & 00:09:58 & 2.9 out of 5\\
		AW6082 & AutoML2 & 0.2247 $\pm$ 0.0871 & 00:09:58 & 2.6 out of 5\\
		AW6082 & SIAMOR  &  0.2497 $\pm$ 0.0606 & 00:04:00 & 0 out of 5\\
		\noalign{\smallskip}\hline\noalign{\smallskip}
		AW7075 & SIASCOR & 0.7965 $\pm$ 0.1710 & 00:00:23 & 0 out of 7\\
		AW7075 & AutoML1 & 0.6071 $\pm$ 0.2291 & 00:09:58 & 4.8 out of 7\\
		AW7075 & AutoML2 & 0.6071 $\pm$ 0.2265 & 00:09:58 & 4.9 out of 7\\
		AW7075 & SIAMOR  & 0.7967 $\pm$ 0.1711 & 00:00:03 & 0 out of 7\\
		\noalign{\smallskip}\hline\noalign{\smallskip}
	\end{tabular}
\end{table}

\begin{figure}[t]
	\sidecaption[t]
	\includegraphics[scale=.375]{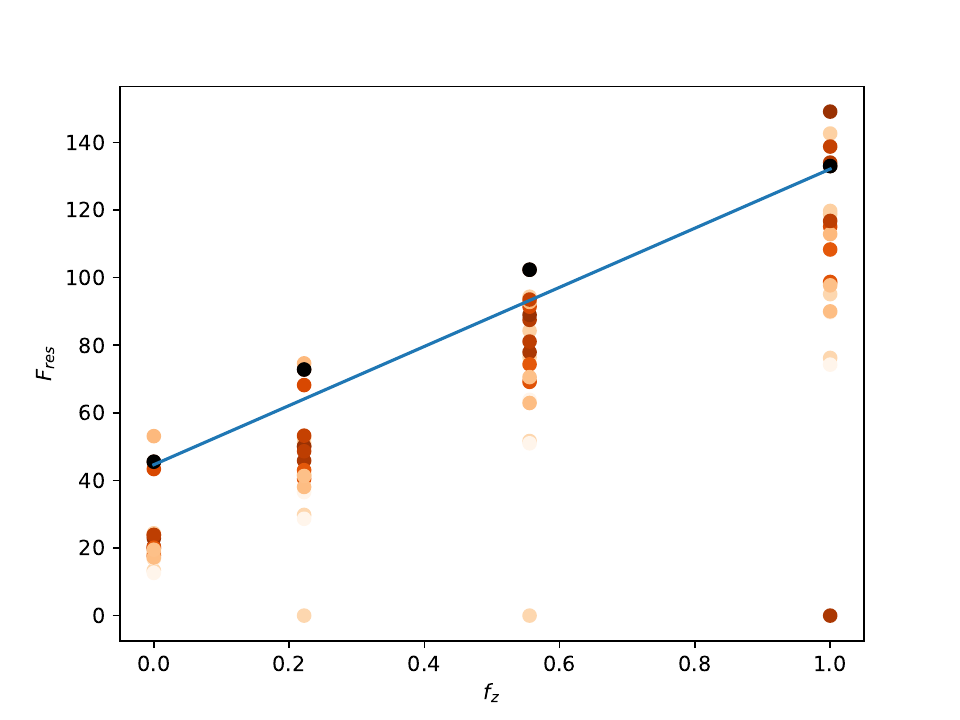}
	\includegraphics[scale=.375]{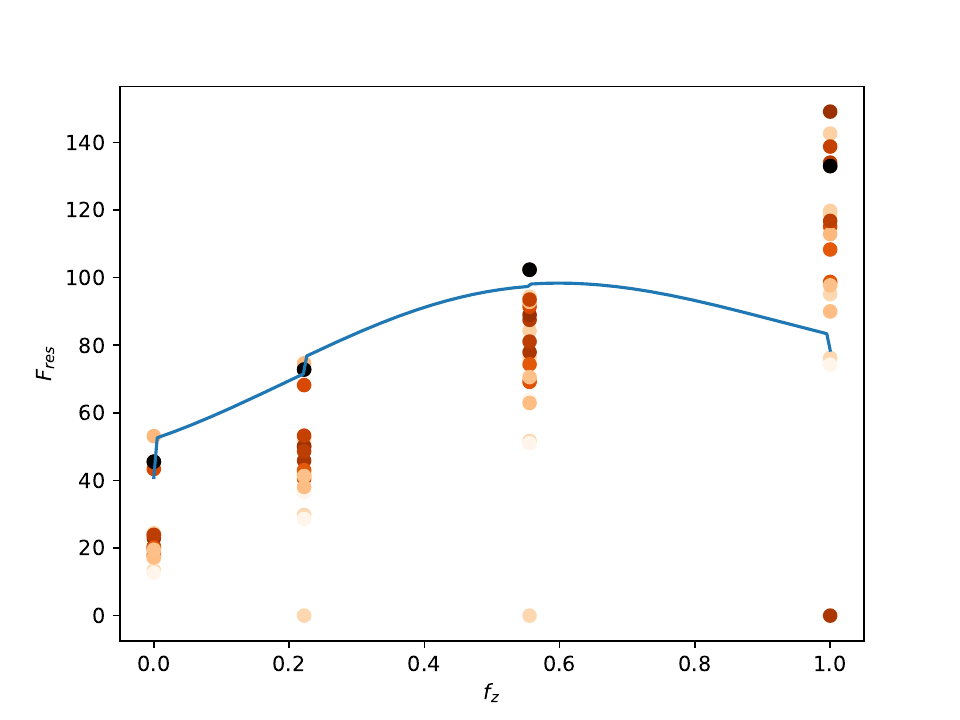}
	\caption{The left  graph depicts the shape of the SIASCOR model and the right graph the shape of the auto-sklearn model both along the $f_z$ variable fixed at the point $x=(0, 0.11, 1)$. The black points represent the data points located along this axis and the orange points are the projected data points onto the $f_z$ dimension. The darker the points are, the closer the Euclidean distance is between the data point and its projection. We can see that the auto-sklearn model strongly violates the monotonicity and convexity. Both models were trained on 90\% of the data that are visible in the graph.}
	\label{fig:violation-milling}
\end{figure}

Table 3 and 4 show the results for all trained models in the milling case. We see again that the two purely data-driven methods violate most of the expected shape behavior. Figure \ref{fig:violation-milling} shows one example of such a shape violation. Here we can see how severe the shape violations can be when shape constraints are not explicitly imposed. More specifically, the AutoML1 for material ``EN-AW7075'' and output  $F_c$ violates both the monotonicity and the convexity constraint along the $f_z$ direction. This time, not only the SIASCOR model but also the SIAMOR model had no violations. Moreover, the training times for both expert-based methods were too low to see a pattern. Both algorithms were fast due to the low dimensionality of $X$ and $W$. Again, the test errors were close to each other considering the ranges of the output.

\section{Artificial Application Example}
\label{sec:artificialexample}

In the previous section, we have considered real-world use cases of small data regression. We have seen that purely data-driven models fail to provide reliable models in terms of shape-compliance. Another indicator for the quality of models is the generalization error, i.e. the  error on data points not contained in the training set. In this section, we examine the generalization error for SIASCOR and compare it to AutoML, SIAMOR, and Ridge regression.

One of the challenges with real small data problems is to compute good estimates of the generalization error with the limited data available. The standard approach to is conduct a cross-validation on the data set, as we have done in the previous real-world use cases. However, this approach limits the focus on the data set and  we cannot infer the performance in regions with few or no data. Additionally, the authors in \cite{MARTENS1998} show that cross-validation can lead to over-optimistic estimates, which may not be indicative of the actual generalization error.

To better understand the generalization error, we construct an artificial application example where we can sample as many data points as we desire. Specifically, we construct an artificial, separable function given by $f^{\text{toy}}:[0, 1]^5 \to \mathbb{R}$,
\begin{align}
	f^{\text{toy}}(x) = \sum_{i=1}^5 f_i^{\text{toy}}(x_i) + 0.1
\end{align}
where 
\begin{equation}
	\begin{aligned}
		f_1^{\text{toy}}(x_1) &= 0.12(x_1 - 0.5)^3, & f_2^{\text{toy}}(x_2)&= 0.002/(x_2 + 0.1)^{2}, \\ 
		f_3^{\text{toy}}(x_3) &= 0.1(x_3 - 0.6)^2(x_3 - 2.4)^2, & f_4^{\text{toy}}(x_4)&= 0.02(x_4 - 0.6)^2(x_3 - 2.4)^2, \\
		f_5^{\text{toy}}(x_5)&= 0.02(x_4 - 1.1)^2(x_3 - 3)^2. & 
	\end{aligned}
\end{equation}
The artificial function is polynomial in all but the second dimension and satisfies all the shape constraints of the brushing example, which can be verified easily.

As in the section before, we compare SIASCOR to SIAMOR and AutoML. The ansatz functions of SIASCOR and SIAMOR are polynomial, while the AutoML conducts a model selection over various ansatz function classes, including polynomial functions.  
In view of the nearly polynomial structure of the artificial function, we compare the shape-constrained models to unconstrained ridge polynomial regression. This way, we can demonstrate that it is the shape constraints that lead to better generalization power and not just the knowledge of the polynomial structure.

To generate the artificial data set, we sample 30 points uniformly from the domain $X = [0,1]^5$, evaluate $f^{\text{toy}}$ at these points, and add Gaussian noise with standard deviation $\sigma = 0.03408$ to the outputs. This corresponds to the expected error in the brushing case. Indeed, experts expect 5-10\% error measuring the roughness. Choosing $\sigma$ as above implies that 95\% of the data have an error below 10\%.

Given the data set and the shape constraints, we formulate the shape-constrained regression problem Example 6.1 in \cite{Schmid.2021}. First, we set the maximal degrees of our anisotropic polynomial as (1, 5, 2, 2, 2). Accordingly, we have $W=[10^{-5}, 10^{5}]^{136}$. Then, we set the regularization parameter as $\lambda =0.01$ for SIASCOR and ridge and $\lambda = 0.05$ for SIAMOR after having conducted a hyperparameter optimization via a ten-fold cross-validation over the dataset. With these settings, the SIP is formulated and we can run SIASCOR and SIAMOR to solve it using the same, remaining settings as in the brushing section. In the case of SIASCOR, we solve the SIP with optimality precision $\delta = 10^{-5}$.

With the hyperparamters settings above, we first conduct another ten-fold cross-validation on the 30 data points. The test errors of the cross-validation are the predictive power estimates that we would have if we had only the data set available. For the cross-validation, we measure the RMSE in every fold and then compute the mean and standard deviation over all folds. Afterwards, we estimate the generalization error. For this, we retrain the models on the entire data set and measure the RMSE on 5000 points sampled from the ground truth function without noise. The results are listed in Table \ref{tab:artificial}.
\begin{table}[!t]
	\caption{Artificial Example: Comparison of SIASCOR, Auto-Sklearn with custom settings, SIAMOR, and Ridge regression. The models were trained on each fold of a 10-fold cross-validation. Then, we measured the average cross-validated (CV) test error. Moreover,  retrained the models on the entire data set and counted the shape violations. Ten shape constraints were considered. Lastly, we measured the generalization error by sampling 5000 points from the ground truth and computing the error to the prediction of the models.}
	\label{tab:artificial}       
	%
	%
	\begin{tabular}{p{3cm}p{2.75cm}p{2.75cm}p{2.75cm}}
		\hline\noalign{\smallskip}
		Model & CV Test Error & Shape Violations & Generalization Error \\
		\noalign{\smallskip}\svhline\noalign{\smallskip}
		SIASCOR & 0.04432 $\pm$ 0.0153 & 0 & 0.032\\
		AutoML2 & 0.05080 $\pm$ 0.0194 & 8 & 0.047\\
		SIAMOR & 0.04230 $\pm$ 0.0153 & 0 & 0.034\\
		Ridge  & 0.05750 $\pm$ 0.0156 & 8 & 0.065\\
		\noalign{\smallskip}\hline\noalign{\smallskip}
	\end{tabular}
\end{table}

First and above all, we observe in Table \ref{tab:artificial} that the shape-constrained models achieve better generalization error than unconstrained models such as AutoML or pure ridge regression. Interestingly, the cross-validated test error on the data set is misleading in the case of SIAMOR because it appears to generalize better than SIASCOR. The generalization error is, however, lower for SIASCOR. This is in line with the observations in \cite{MARTENS1998} where the authors have shown that cross-validation usually slightly underestimates the true generalization error. Nonetheless, for the unconstrained models the cross-validation error is indicative of the inferior generalization error. 

In conclusion, the use of shape constraints is shown to lead to better generalization power in the artificial small data regression problem that we considered. This can be inferred from the lower generalization errors of SIASCOR and SIAMOR over AutoML. By considering Ridge regression that has the same ansatz functions, we conclude that it is indeed the shape constraints that lead to a better performance and not the choice of ansatz functions. All in all, the evaluation of the generalization error on this artificial data set emphasizes another benefit of enforcing shape constraints


\section{Conclusion}
\label{sec:conclusion}
Reliable models are important in situations with insufficient data. Take for instance the quality prediction in manufacturing. Data generation is expensive because it involves costly experiments. Therefore, data is scarce and, on top of that, noisy due to measurement errors. During adaptive process control, quality prediction is crucial to avoid wrong decision-making that comes along with high costs. Therefore, practitioners rely rather on their prior knowledge than on models trained solely on the available data. Informed learning aims to get the best from both sides: inferring quantitative information from data while being in accordance with prior knowledge. 

In this work, we incorporated shape knowledge into the training of machine learning models. First, we summarized the SIASCOR method from \cite{Link.2022} by formulating the shape-constrained regression problem and discussing SIP algorithms for solving it. We also recalled the ISI methodology as a general method to inspect models for shape compliance or non-compliance, to elaborate and specify shape constraints, and to incorporate these shape constraints into the chosen regression model. We point out that SIASCOR is only one possible way of incorporating shape constraints. In principle, other methods of integrating shape knowledge can be used within the ISI methodology.  Then, we considered three application examples from manufacturing: brushing, press hardening, and milling. The latter has not been applied to shape-constrained regression yet, so we used ISI to formalize shape constraints. Then, we applied the SIASCOR method with the simultaneous algorithms from \cite{Schmid.2021} to obtain shape-compliant models that fit the data. The SIASCOR method was then compared to two purely data-driven AutoML approaches and another shape-constrained regression method SIAMOR that solves the corresponding SIP differently.  During a 10-fold cross-validation, we created ten different model functions for each approach and compared their shape compliance (or more precisely the average number of shape constraints they violated), their average training time, and their average test error. 

In the comparative study with real-world application examples, we have seen that, in general, purely data-driven models trained on a few data do not infer shape knowledge just from the data. Even methods like auto-sklearn that choose the best model out of many do not seem to perform satisfactorily outside the data set. In fact, the resulting models behave contrary to physical laws repeatedly.  Consequently, we do not recommend using purely data-driven methods to create reliable models when the available data is limited. Furthermore, just by looking at the test errors computed with cross-validation, one is tempted to think that AutoML models generalize slightly better. But as we have seen, low test errors are misleading because the final models violate the shape constraints. In contrast, the final SIASCOR model was, as expected, shape-compliant.

In the comparative study with the artificial example, we could also show that the use of shape constraint lead to more powerful models regarding generalization error. When trained on a small data set, purely data-driven methods, as AutoML or Ridge regression, do not perform well outside the data set.

In conclusion, we recommend leveraging shape knowledge when dealing with insufficient data.  However, not all shape-knowledge-based approaches guarantee definite shape compliance of the final model. For instance, SIAMOR fails to rigorously impose shape constraints. From a theoretical point of view, this is to be expected but we have also seen it in practice. Besides, there are no mathematical guarantees for the termination of the algorithm and inferior properties for optimality. In contrast to SIAMOR, SIASCOR is a theoretically sound and practically reliable method to impose hard shape constraints as we have seen in the experiments. Aside from the shape compliance, we observed in our applications that for a higher number of inputs and parameters the SIAMOR is slightly faster than SIASCOR. However, the discrepancy between the training times was not large enough for proper assessment. Even if there was a significant gain in computation time, it would only justify to use SIAMOR for exploring shape-constrained regression but not for imposing shape constraints on the final model.

A follow-up study should reconsider the adaptive feasible-point algorithm that solves the SIP to improve computation time. One theoretical aspect is to analyze the time complexity of the algorithms. This is crucial for better scalability with respect to an increasing number of model parameters and input variables. Another, more applied aspect is, of course, an efficient implementation. It is also interesting to see how the generalization error of  SIASCOR changes with respect to different noise levels and to different numbers of data. This analysis would help to better estimate when to use the SIASCOR. An additional line of research can be to develop an elaborated way to find all regions where shape constraints are violated. This way, we cannot only visualize all shape violations but, on top of that, accelerate the process of finding the right shape constraints with ISI. One could also, in a comparative study, analyze if, and to which extent, SIASCOR improves extrapolation. Going in a similar direction, a tool that helps to assess whether some enforce shape constraints are in too much conflict with the data and should therefore be excluded would be very handy as we have seen in the milling example. One further idea is to find a way to compute the confidence intervals of the SIASCOR model to spot regions of high and low variance.

%
%
\bibliographystyle{acm}
\bibliography{bibliography.bib}

\end{document}